\documentclass[11pt, a4paper]{article}
\usepackage{amsmath, amsthm, amssymb, eucal, setspace}

\makeatletter
\renewcommand\section{\@startsection{section}{1}{\z@}%
 						{-3.5ex \@plus -1ex \@minus -.2ex}% negative means
						{2ex \@plus.2ex}% 		positive means vertical skip
						{\large\bfseries}}
\renewcommand\subsection{\@ifstar
						{\setcounter{subsection}{\value{equation}}
					\@startsection{subsection}{2}{\z@}
                          {1.75ex \@plus.5ex \@minus.2ex}%
                           {-.4em}		% negative means horizontal
					\textit*}
					{\setcounter{subsection}{\value{equation}}
						\stepcounter{equation}
					\@startsection{subsection}{2}{\z@}
                          {1.75ex \@plus.5ex \@minus.2ex}%
                           {-.4em}		% negative means horizontal
					\textit}}
\def\@seccntformat#1{\@ifundefined{#1@cntformat}%
	{\csname the#1\endcsname\quad} 
	{\csname #1@cntformat\endcsname}} 
\def\section@cntformat{\thesection.~} 
\def\subsection@cntformat{(\thesubsection)\ }
\renewcommand*\l@section{\mdseries\small\@dottedtocline{1}{1.5em}{2em}}
\makeatother

\numberwithin{equation}{section}
\theoremstyle{plain}
\newtheorem{maintheorem}{Theorem}
\swapnumbers
\newtheorem{theorem}[equation]{Theorem}
\newtheorem{corollary}[equation]{Corollary}
\newtheorem{lemma}[equation]{Lemma}
\newtheorem{proposition}[equation]{Proposition}
\theoremstyle{definition}

\theoremstyle{remark}

\newcommand{\cE}{\mathcal{E}}
\newcommand{\cO}{\mathcal{O}}
\newcommand{\frg}{\mathfrak{g}}
\newcommand{\frt}{\mathfrak{t}}
\newcommand{\frb}{\mathfrak{b}}
\newcommand{\frz}{\mathfrak{z}}

\newcommand{\bQ}{\mathbb{Q}}
\newcommand{\bZ}{\mathbb{Z}}
\newcommand{\Ind}{\mathrm{Ind}}
\newcommand{\Id}{\mathrm{Id}}
\newcommand{\Tr}{\mathrm{Tr}}

\newcommand{\GL}{\mathrm{GL}}
\newcommand{\SL}{\mathrm{SL}}
\newcommand{\frM}{\mathfrak{M}}
\newcommand{\mi}{\mathrm{i}}
										
\begin{document}

\title{\textbf{The Newstead-Ramanan conjecture for Chern classes}}
\author{Constantin Teleman\footnote{Partially supported by EPSRC grant 
S06165/01.}
 \and Christopher T.~Woodward\footnote{Partially supported by NSF grant 
 DMS/0093647.}}
\date{ 20 December 2005}
\maketitle
\centerline{\emph{Dedicated to the memory of Raoul Bott}}

\begin{quote}
\abstract{\noindent Newstead and Ramanan conjectured the vanishing of the 
top $2g-1$ Chern classes of the moduli space of stable, odd degree vector 
bundles of rank $2$ on a Riemann surface of genus $g$. This was proved 
by Gieseker \cite{gies}, while an analogue in rank $3$ was settled by 
Kiem and Li \cite{kiemli}. We generalise this to the vanishing of the 
top $(g-1)\ell$ rational Chern classes of the moduli space $M$ of 
stable principal bundles with semi-simple structure group of rank 
$\ell$, whenever $M$ is a compact orbifold.}
\end{quote}

%%%%%%%%%%%%%%%%%%%%%%%%%%%%%%% Introduction %%%%%%%%%%%%%%%%%%%%%%%%
\section*{Introduction}
Let $G$ be a complex reductive group and $M$ the moduli space of stable 
principal $G$-bundles on a compact Riemann surface of genus $g$. While 
$M$ is smooth (in the orbifold sense), it is usually not compact, but 
is compactified by adding the semi-stable bundles. In important special 
cases, semi-stable bundles are stable and then $M$ is a compact orbifold. 
This happens when $G=\GL(n)$ for the components of degree co-prime to 
$n$, or else if we enrich the bundle with a sufficiently generic 
parabolic structure at a marked point on the surface.\footnote{These 
examples are unified by noting that the moduli of stable vector bundles 
degree $d$ is also the space of stable vector bundles of degree $0$ but 
with parabolic structure defined by the conjugacy class $\exp(2\pi\mi 
d/n)\cdot\Id$ in $\GL(n)$.} 
Henceforth, we place ourselves in one of these favourable situations. 
Let $\ell^{ss}$ and $\ell^c$ be the semi-simple and central ranks of $G$. 
We prove here the following generalisation of an old conjecture of 
Newstead and Ramanan \cite{new, ram}.

\begin{maintheorem}
The top $(g-1)\ell^{ss} + g\ell^c $ rational Chern classes 
of $M$ vanish. 
\end{maintheorem}

\noindent For rational cohomology, we can pass to finite covers with 
impunity \cite[\S7]{ab} and split $G$ as a torus of rank $\ell^c$ and 
simple groups; so the only content of the theorem concerns $\ell^{ss}$. 
We actually prove an equivalent result in topological $K$-theory. Let 
$G$ be semi-simple of rank $\ell$.

\begin{maintheorem}
The top $(g-1)\ell$ rational Grothendieck $\gamma$-classes of $M$ vanish. 
\end{maintheorem}
  
\noindent The Grothendieck $\gamma$-classes are the $K$-theoretic Chern 
classes; we recall them below, along with the equivalence of Theorems 
1 and 2. In some cases, such as for $G=\GL(n)$ or $\SL(n)$, we known 
thanks to Atiyah and Bott \cite{ab} that $M$ is free of homology (and 
therefore $K$-theory) torsion, and we get an integral result. Experience 
with equivariant $K$-theory suggests that $K(M)$ might be torsion-free 
whenever $\pi_1G$ is so, but this does not seem to be known.

To prove Theorem 2, we pair the total $\gamma$-class $\sum t^k\gamma^k$ 
of $T^*M$ against any test class in $K^0(M)$ and show that we obtain a 
polynomial of degree no more than $\dim M -(g-1)\ell$. The computation 
is straightforward, but relies on two substantial results. The first 
is the \textit{index formula} of \cite{tw}, an explicit Verlinde-like 
formula for the index (direct image) of ``almost any" $K$-theory class 
over the stack $\frM$ of $G$-bundles. The various parabolic structures 
on bundles leads to numerous minor variations in the formula; while it 
takes no effort to write down any of them, a shortcut allows us to 
restrict attention to the \emph{full-flag} structures (Corollary~\ref
{borelenough}). The second ingredient relates the indexes over $\frM$ 
and over the moduli space $M$ of vector bundles twisted by 
K\"ahler differentials $\cite[\S7]{tel}$. 

Such a proof by Poincar\'e duality is not new, cf.~Zagier \cite{zag} 
for $\SL(2)$, but the cohomological formulae for integration over 
$M$ turned out to be unwieldy. One way to understand why our index 
formula works well is to interpret it as a localisation to the stack 
of principal bundles under the maximal torus $T\subset G$ \cite
[\S2.17]{tw}; there, the tangent complex to $\frM$ has a trivial 
summand of rank as predicted by the vanishing. Still, the method 
has its limits. Thus, Theorem~2 leads one to ask whether the 
$\gamma$-classes vanish in \textit{algebraic} $K$-theory. The 
absence of an obvious counter-candidate leads one to suspect so, 
but the current method fails.

%%%%%%%%%%%%%%%%%%%%%%%%%%%%%% Section 1 %%%%%%%%%%%%%%%%%%%%%%%%
\section{The $\gamma$-classes}
For a complex vector bundle $V$ of rank $r$ over a compact space $X$, 
let $\lambda^k(V)\in K^0(X)$ be its $k$th exterior power and define 
the following polynomials of degree $r$ with coefficients in $K^0(X)$:
\[\begin{split}
\lambda_t(V) =& \sum\nolimits_k t^k\lambda^k(V),\\ 
\gamma_t(V) =& \sum\nolimits_k t^k \gamma^k(V) := (1-t)^r \lambda_
{t/(1-t)}(V).
\end{split}\]
The second relation defines the classes $\gamma^k(V)\in K^0(X)$. Note 
that $\gamma_t(V\oplus W) =\gamma_t(V)\cdot\gamma_t(W)$ for vector bundles 
$V$ and $W$, while $\gamma_t(L) = (1-t) + tL$ for a line bundle $L$; by 
the splitting principle, these conditions determine $\gamma_t$.  Also, 
$\gamma^1(L) = L -1$, the $K$-theory Euler class of the line bundle, and 
in this sense $\gamma_t$ is the total $K$-theory Chern class. The next 
exercise is included for the reader's convenience.

\begin{proposition}
The following assertions about the vector bundle $V$ are equivalent.
\begin{trivlist}\itemsep0ex
\item (i) The top $d$ rational Chern classes vanish.
\item (ii) The top $d$ rational $\gamma$-classes vanish.
\item (iii) The polynomial $\lambda_t(V)\in K^0(X;\bQ)[t]$ is divisible 
by $(1+t)^d$.
\end{trivlist}
\end{proposition}

\noindent When $K^0(X;\bQ)$ satisfies Poincar\'e duality with respect to a 
map $\Ind: K^0(X) \to \bQ$, these conditions are equivalent to

\noindent (iv) For every $W\in K^0(X)$, $\Ind(\lambda_t(V)\cdot W)\in 
\bQ[t]$ vanishes to order $d$ or more at $t= -1$. 

\begin{proof}
Equivalence of (ii) and (iii) is clear from the inversion formula 
$\lambda_t = (1+t)^r\gamma_{t/(1+t)}$. Observe now that in the ring 
$R$ of symmetric power series in variables $x_1,\dots,x_r$, the 
ideal $(e_{r-d+1}, \dots, e_r)$ generated by the top $d$ elementary 
symmetric functions is the intersection of $R$ with the ideal 
$(x_{r-d+1},\dots,x_r)\in \bQ[[x_1,\dots,x_r]]$. The transformation 
$x_k\mapsto y_k:= e^{x_k} -1$ defines an automorphism of $\bQ[[x_1,
\dots,x_r]]$ which preserves $(x_{r-d+1},\dots,x_r)$. It follows 
that $(e_{r-d+1}, \dots, e_r)$ agrees with the ideal of the top $d$ 
elementary symmetric functions in the $y_k$. We now let $x_k$ be the 
Chern roots of $V$; then, $ch\,\gamma_t(E) = \prod(1 + ty_k)$, so the 
$\gamma$-classes are the elementary symmetric functions in the $y$'s, 
and we conclude that $\text{(i)}\Leftrightarrow\text{(ii)}$. 
\end{proof}

%%%%%%%%%%%%%%%%%%%%%%%%%%%%%% Section 2 %%%%%%%%%%%%%%%%%%%%%%%%
\section{$K$-theory of $\frM$ and $M$}

Let $\frM$ be the stack of all algebraic $G$-bundles over $\Sigma$. 
Its homotopy type is that of the space of continuous maps from $\Sigma$ 
to $BG$. For equivariant $K$-theory, it is more natural to assign 
to $\frM$ an equivariant homotopy type, given by the conjugation 
action of the maximal compact subgroup $G_K\subset G$ on the space of 
\emph{based} continuous maps from $\Sigma$ to $BG$. The last space is 
a principal fibration, with fibre the space of based loops on $G$, 
over a product of copies of $G$ \cite{ab}. By $K^*(\frM)$ we will 
mean the $G_K$-equivariant $K$-theory of the space of based maps. 

Even though $\frM$ has an underlying algebraic structure, one cannot 
quite define an ``index map" from $K^*(\frM)$ to $\bZ$, because 
$\frM$ has infinite type. This is most obvious for $G= \GL(1)$, when 
there are infinitely many components. Nonetheless, in \cite{tw} we 
defined and computed an index map on a sub-ring of $K^*(\frM)$. In 
this section, we recall the generating classes and the index formula, 
adapted to include the \emph{K\"ahler differentials} of $\frM$. 
We remark that the fine points of the definition of the index need 
not concern us here, because the stack $\frM$ is only an intermediate 
step: what we need is the index over $M$, where the definition is clear.

\subsection{$K$-theory generators.}
The universal principal $G$-bundle over $\Sigma\times\frM$ associates 
to a $G$-representation $V$ a vector bundle $E^*V$. Slant product with 
cycles $C\in K_{0,1,2}(\Sigma)$ gives classes $E^*_CV \in K^{0,-1,-2}
(\frM)$. For the degree $0$ class, we choose a point $x$, while in 
degree $2$ we will use the fundamental class of $\Sigma$; if $\pi$ 
is the projection along $\Sigma$ and $K$ the relative canonical bundle, 
$E^*_\Sigma V$ is the class of $R\pi_*(E^*V\otimes\sqrt K)$ in $K^0\cong 
K^{-2}$ (by Bott periodicity). 

If $G$ is simply connected, one can show by using Chern characters and 
the completion theorem that the $E^*_CV$ generate a dense sub-ring of 
 $K^*(\frM;\bQ)$. The Atiyah-Bott surjectivity argument \cite{ab} also 
shows that we generate $K^*(M;\bQ)$ this way, provided that all semi-stable 
bundles are stable. When our bundles are enhanced with a complete flag 
at $x$ (see \S\ref{parab} below; this is the only parabolic type we 
will need), we can define $E^*_xV$ from any representation $V$ of the 
\textit{Borel subgroup} of $G$, and these again generate the rational 
$K$-theory of the space of stable bundles. Naturally, for the latter, 
we only need to consider completely reducible representations $V$. 

\subsection{The stack and the space.}  Denote by $\Omega^\bullet$ the 
K\"ahler differentials over the stack $\frM$ and $\Omega^\bullet(h)$ 
their twist by the line bundle $\cO(h)$ with Chern class $h\in H^2
(\frM;\bZ)$ ($=\bZ$, when $G$ is simple and simply connected). Let also 
$\cE$ be a monomial $\bigotimes_i E^*_{C_i}V_i$. The bundles $\Omega^\bullet
(h)\otimes \cE$ might fail to descend to the moduli space $M$ (there can 
be obstructions from the finite stabilisers); if so, what is implicitly 
meant below are their invariant direct images, from the sub-stack of stable 
bundles to $M$. The space and the stack are related by the following 
proposition, in which we may allow matching parabolic structures on $M$ 
and $\frM$. 

\begin{proposition}\label{space}
When all semi-stable bundles are stable, the cohomologies of $\Omega^
\bullet(h)\otimes \cE$ over the stack $\frM$ and the space $M$ agree 
for sufficiently large $h$,\footnote{If $\lambda$ is the sum of the 
highest weights of the $V_i$ and $\alpha_0$ the highest root of $\frg$, 
any $h>b(\lambda,\alpha_0)$ in the basic bilinear form $b$ on $\frg$ 
will do, save that a small adjustment by the parabolic weights may be 
needed.} independent of the degree of the differentials. 
\end{proposition}

\noindent The co-factor $\cE(h)$ of $\Omega^\bullet$ will be a ``test 
class" in $K^0$. In the next section, we will check that the index of 
the formal series $\sum t^k\Omega^k\otimes\cE(h)$ over $\frM$ vanishes 
to order $(g-1)\ell$ at $t=-1$. Since the index over the space $M$ is 
quasi-polynomial in $h$, this behaviour applies to all $h$, if it does 
to large $h$. As our test classes $\cE$ span $K^*(M)$, we will have 
proved Theorem~2.

\begin{proof}
We use the local cohomology vanishing results \cite[\S7]{tel} along 
unstable strata of $\frM$. Specifically, consider the projection 
$\Sigma^m\times \frM \twoheadrightarrow \Sigma^m$, with the fibre-wise 
stratification of $\frM$. We apply the vanishing results to the direct 
image of the bundle $\Omega^\bullet(h) \otimes (\boxtimes_i E^* V_i)
\otimes E^*_xU$ along $\frM$, where the $V_i$ are the $G$-representations 
defining the factors of $\cE$. The algebraic isomorphism type of the 
vector bundles $E^*V_i$ along the fibres $\frM$ varies, but their weights 
along unstable strata do not, and they alone enforce the vanishing. So 
the direct images along $\frM$ and $M$ agree. The resulting class in 
$K^0(\Sigma^m)$ can now be integrated along the product of the cycles 
$C_i$ and shows the equality of the indexes of $\Omega^\bullet(h)\otimes 
\cE$ over $\frM$ and $M$.    
\end{proof}

\subsection{Parabolic structures.}\label{parab}
For some background on moduli spaces and stacks of bundles with parabolic 
structures, see \cite[\S9]{tel}. Our main example will be the stack 
$\frM(B)$ of bundles with a ``full flag" (meaning a coset for the Borel 
subgroup $B\subset G$) in the fibre over $x\in\Sigma$. This is a 
$G/B$-bundle over $\frM$, associated to the universal $G$-bundle
 restricted to $\{x\}\times\frM$.

Any other stack $\frM'$ of bundles with a parabolic structure at $x$ is 
the base of a fibre bundle $\pi:\frM(B)\to \frM'$. The fibres are 
homogeneous spaces for a Levi subgroup of the loop group of $G$. Since 
$R\pi_*\pi^* = \Id$ on $\cO$-modules, the index of a coherent sheaf 
over $\frM'$ equals that of its lift to $\frM(B)$. It suffices therefore, 
in principle, to understand the index formula over $\frM(B)$; this 
will be given in the next section. However, our interest lies mainly 
with the K\"ahler differentials, and their effect on the push-forward 
can be described explicitly. 

Over $\frM(B)$, we have a distinguished triangle of tangent complexes 
\[
T_\pi\frM(B) \to T\frM(B) \to \pi^*T\frM' \to T_\pi\frM(B)[1].
\]
When $\frM' =\frM$, the tangent bundle $T_\pi\frM(B)$ along the fibres 
is the evaluation bundle $E_x^*(\frg/\frb)$. Splitting the sequence 
gives an equality in $K$-theory, 
\[
\lambda_t\left((T^*\frM(B)\right) = \lambda_t\left(T^*_\pi\frM(B)\right)
	\otimes \pi^*\lambda_t(T^*\frM').
\]
The fibres of $\pi$ are smooth and proper, with Hodge type $(p,p)$. Hodge decomposition gives 
\[
R\pi_*\left[\lambda_t\left(T^*\frM(B)\right)\otimes\pi^*\cE(h)\right] = 
	\lambda_t(T^*\frM')\otimes\cE(h)\cdot \sum (-t)^p b_{2p}
\]
where the $b_{2p}$ are the Betti numbers of the fibre. For $t=-1$, the 
last factor becomes the Euler characteristic, which is positive and hence
and does not affect the vanishing order of the index. 
%We summarise this insight in the following 

\begin{corollary}\label{borelenough}
If Theorem 2 holds for moduli of bundles with full-flag structures, 
then it holds for all parabolic structures. \qed
\end{corollary}

%%%%%%%%%%%%%%%%%%%%%%%%%%%%%%%%%%%% Section 3 %%%%%%%%%%%%%%%%%%%%%%%%%
\section{The index formula} 

We now describe the ingredients of our index formula: a function $\theta$ 
on the maximal torus $T\subset G$ and a subset $F$ of points in its domain. 
Note that the multiple $hb$ of the basic bilinear form $b$ defines a 
homomorphism $e^{hb}:T\to T^*$ to the dual torus. Next, the differential 
$d\Tr_V$ is a regular $\frt^*$-valued function on $T$, and so $e^{s
\cdot d\Tr_V(.)}$ is a map from $T\to T^*$, depending formally on $s$. 
Similarly, $1+te^\alpha$ is a function on $T$, so $(1+te^\alpha)^
\alpha$ is a $T^*$-valued map. Set 
\begin{equation}\label{chi}
\chi= e^{(h+c)b + s\cdot d\,\Tr_V(.)}\cdot 
	\prod_{\alpha>0}\left[\frac{1+te^\alpha}{1+te^{-\alpha}}\right]^\alpha
: T\to T^* 
\end{equation} 
($c$ is the dual Coxeter number of $\frg$). The set $F$ of solutions 
of the equation
\begin{equation}\label{points}
\chi(f) = 1 \in T^*.
\end{equation}
depends on $h,V$ and on the formal variables $s,t$. We denote by $F^{reg}$ 
the subset of solutions which are \textit{regular} (as $G$-conjugacy 
classes) at $s=t=0$.

Call $H(f)$ the differential of $\chi$ at $f\in T$; the notation $H$ 
stems from its agreement with the Hessian of the function on $\frt$
\[
\xi \mapsto \frac{h+c}{2}\,b(\xi,\xi) + s\Tr_V(e^\xi) - 
\Tr_\frg\left(\mathrm{Li}_2(te^\xi)\right), 
\]
where $\mathrm{Li}_2$ is Euler's dilogarithm. Using the metric $(h+c)b$ 
to identify the Lie algebra $\frt$ of $T$ with its dual, we can convert 
$H$ to a matrix and define
\begin{equation}\label{theta}
\theta(f)^{-1} = |F|\cdot\prod_\alpha 
		\frac{1+te^\alpha}{1-e^\alpha}\cdot\det H(f). 
\end{equation}
Note that $\det H =1$ at $s=t=0$.

\begin{theorem}\label{index} For $h\ge 0$, we have the index formula
\[
\sum\nolimits_{k\ge 0} t^k\cdot \Ind\left(\frM; \Omega^k(h)\otimes
		\exp[s E^*_\Sigma V]\otimes E^*_x U\right) = 
		(1+t)^{(g-1)\ell}\sum\nolimits_f \theta(f)^{1-g}\cdot \Tr_U(f),
\]
with $f\in F^{reg}$ ranging over a complete set of Weyl orbit 
representatives. 
\end{theorem}

\subsection{Odd generators.} \label{skew}
Each quadratic form $H(f)$ also gives a bilinear form $\langle.|.\rangle
_f$ on the dual $\frt^*$ of the Lie algebra. If $C_{1,2}$ are two 
$1$-cycles on $\Sigma$ with intersection product $\#C_1 C_2$, then an 
additional quadratic monomial factor $E^*_{C_1}W_1\wedge E^*_{C_2}W_2$ 
in the coefficients of \eqref{index} contributes a co-factor of $\#C_1C_2 
\cdot \langle d\Tr_{W_1}(f) | d\Tr_{W_2}(f) \rangle_f$ to the $f$-term in 
the right-hand sum. This recipe defines a skew contraction procedure on 
odd-generator monomials, which captures the answer in general: we sum, 
over all possible complete contractions, the product of the corresponding 
pairings, and insert the result as a co-factor for each term \cite
[Thm.~2.15]{tw}.

\begin{proof}[Proof of the index formula.]
The series $\bigoplus_k t^k\Omega^k$ of complexes of K\"ahler differentials 
on $\frM$ is the total $\lambda_t$-class of the complex $R\pi_*(E^*\frg
\otimes K)$, with the adjoint representation $\frg$. The last direct 
image is $E^*_\Sigma\frg\oplus E^*_x \frg^{\oplus(g-1)}$ in topological 
$K$-theory (mind the extra $\sqrt K$ in our bundle). 

We express $\lambda_t$ in terms of Adams operations $\psi^k$,
\[
\lambda_t(\frg) = \exp\left[-\textstyle\sum_{k>0} (-t)^k\psi^k(\frg)/k\right],
\]
and note the relation 
\[
\psi^k (E^*_\Sigma\frg) = \frac{1}{k} E^*_\Sigma\,\psi^k(\frg).
\]
Theorem \ref{index} then refers to the index over $\frM$ of
\[
\cO(h)\otimes \exp\left[sE^*_\Sigma V - 
	\sum\nolimits_{k>0} \frac{(-t)^k}{k^2} 
		E^*_\Sigma\,\psi^k(\frg)\right]
		\otimes E^*_x\lambda_t(\frg)^{\otimes(g-1)}
		\otimes E^*_xU,
\]
which has the type studied in \cite[Thm.~2.11]{tw}. There, we wanted
the solutions $f$ to  
\[
\exp\left[ (h+c)b + s\cdot d\Tr_V - \sum\nolimits_{k>0;\alpha} 
	\frac{(-t)^k}{k} e^{k\alpha}\cdot\alpha\right] = 1;
\]
but this is precisely \eqref{points}. To relate formula \eqref{index} 
with the index formula in \cite{tw}, just observe the pre-factor 
$(1+t)^\ell$ and the factors $1+te^\alpha$ in $\theta^{-1}$ come from 
splitting the character of $\lambda_t(\frg)$ as $(1+t)^\ell\cdot\prod_
\alpha (1+t e^\alpha)$. The enhancement to odd generators is 
\cite[Theorem 2.15]{tw}.   
\end{proof}

\subsection{Variation: full-flag parabolic structures.}
The formula the stack $\frM(B)$ can be derived from Theorem~\ref{index} 
and from the Weyl character formula for the index along the fibres $G/B$ 
of the projection $\frM(B)\to\frM$. The outcome is to replace $\Tr_U$ 
in the right-hand side of \eqref{index} by $\Tr_U\cdot\prod_{\alpha>0}
(1-e^\alpha)^{-1}$, and to sum over all points of $F^{reg}$ instead 
of Weyl orbits. We must also take care to use the appropriate 
differentials; the fibres $G/B$ contribute a factor $\prod_{\alpha>0}
(1+te^\alpha)$.

%%%%%%%%%%%%%%%%%%%%%%%%%%%%%%% Section 4 %%%%%%%%%%%%%%%%%%%%%%%%%%
\section{The limit $t \to -1$}

The factor $(1+t)^{(g-1)\ell}$ already appears in \eqref{index}, so 
we must check that no singularities in $\theta(f)^{1-g}$ cancel the 
vanishing at $t=-1$. To do so, we must study the roots $f_t$ of 
\eqref{points}.

When $h>0$ and $t=s=0$, $\chi$ is an isogeny, so all solutions to 
\eqref{points} are simple. The following lemma ensures that they remain 
simple for all $t\in (-1,0]$ and small $s$. 

\begin{lemma}\label{hessian}
If $h>c$, $s=0$ and $t\in (-1,0]$, the differential $H = d\chi$ is 
non-degenerate on $T_K$.
\end{lemma}

\begin{proof}%[Proof of (\ref{hessian}).] 
With $H_V(f)$ denoting the Hessian of $\Tr_V$ at $f$, we have
\begin{equation}\label{log}
H = (h+c)b + s H_V(f) + t \sum\nolimits_\alpha 
			\frac{e^\alpha}{1+te^\alpha}(f)\cdot\alpha^{\otimes 2}.
\end{equation}
Note that $\alpha^{\otimes 2}$ is negative semi-definite, $t\le 0$ and 
$\Re\frac{e^\alpha}{1+te^\alpha} \ge -1$ for $|e^\alpha| =1$. As 
$\sum_\alpha \alpha^{\otimes 2} = -2cb$, $H$ is bounded below by 
$(h-c)b + sH_V$.
\end{proof}

\noindent Skew-adjointness of $\chi$ for $s=0$ then implies that $F$ stays 
in the \textit{compact} part $T_K$ of the maximal torus for small variations 
in the real time $t$, and thus for all times $t\in [-1,0]$. Non-degeneracy 
of $H$ also shows that the $s$-dependence in \eqref{points} can be solved 
order-by order, when $t\in (-1,0]$. This keeps $F$ in a formal neighbourhood 
of $T_K$. Below, we will find that $H$ remains regular at $t=-1$, so the 
solution can be perturbed analytically with $s$ even at that value of $t$.  

Since no solutions escape to $\infty$, singularities in the $\theta(f)^
{-1}$ can only stem from zeroes in the denominators in \eqref{theta} 
and singularities in $H$. It turns out that certain points of $F^{reg}$ 
\emph{do} wander into the singular locus of $T$ as $t\to -1$, so we 
need to control this.

\begin{lemma}\label{converge}
Let $f_t\in F$ be regular at $t=0$ but singular at $t=-1$. For 
small $x = \sqrt{t+1}$, $f_t$ has a convergent expansion 
\[
f_t = f_{-1}\cdot\exp\left[\textstyle{\sum}_{k>0}\; x^k\xi_k\right].
\] 
Moreover, $\beta(\xi_1)\neq 0$ for any root $\beta$ such that 
$e^\beta(f_{-1}) = 1$.
\end{lemma}
\noindent Thus, the tangent to $f_t$ at $f_{-1}$ is regular in 
the Lie algebra centraliser $\frz$ of $f_{-1}$. We obtain
\begin{equation}\label{limits}\begin{split}
\lim_{t\to -1}\frac{1+te^\alpha}{1-e^\alpha} &= 
	1 \qquad\text{for all roots } \alpha,\\
\lim_{t\to -1}\left(\frac{e^\beta}{1+te^\beta} + 
	\frac{e^{-\beta}}{1+te^{-\beta}}\right) &= -1 - \frac{2}{\beta(\xi_1)^2}
	\quad\text{for roots $\beta$ of }\frz.
\end{split}\end{equation}
In particular, the limiting value $H(f_{-1})$ in \eqref{log} is the positive 
definite form $hb + sH_V + \sum_\beta \frac{\beta^{\otimes 2}}{\beta
(\xi_1)^2}$, summing over roots of $\frz$. This excludes any 
singularities in \S\ref{skew}.

\begin{proof}[Proof of (\ref{converge}).]
At $t=-1$, equation \ref{points} can be simplified\footnote{We are using 
the relations $(-1)^{2\rho}=1$ for simply connected $G$, as well 
as $\sum_{\alpha>0}\alpha^{\otimes 2} = -c$.} to 
\begin{equation}\label{at-1}
\exp[hb + s\cdot d\Tr_V(.)]= 1;
\end{equation}
however, the cancellation involved conceals multiple solutions on the 
singular locus in $T$. The latter partitions $T_K$ into alcoves that are 
simply permuted by the Weyl group. We claim that 
\begin{itemize}\itemsep0ex
\item Each \emph{singular} solution of \eqref{at-1} is a limit of at 
least one solution behaving as in Lemma \ref{converge}.
\end{itemize} 
By Weyl symmetry, there must be such a solution from each adjacent alcove.
Now, every \emph{regular} solution of \eqref{at-1} is also the limit of 
a regular solution of \eqref{points}: this is because it is the limit of 
\emph{some} solution, and Weyl symmetry ensures that the points of $F$ 
that are singular at $t=0$ stay so at small times, and then at all times. 
Finally, recall that the solutions of \eqref{at-1} in a \emph{closed} Weyl 
alcove are in bijection with the regular solutions of \eqref{points} in 
that alcove. Our claim then accounts for the $t=-1$ limits of \emph{all} 
remaining points of $F^{reg}$ and proves Lemma~\ref{converge}.

To prove the claim, it suffices to find a \emph{formal} solution $f_t$ as 
in the Lemma. As $t$ converges faster than $f_t$ approaches the singular 
locus, the function $\chi$ becomes \eqref{at-1} in the limit, so the 
equation is verified to zeroth order precisely when $f_{-1}$ solves 
\eqref{at-1}. To obtain the constraint on $\xi_1$, 
we differentiate the log of \eqref{chi} in $x$:
\[
\iota(\xi')\left[(h+c)b + sH_V(f_t)\right] +
	\sum_\alpha \frac{[2x + t\alpha(\xi')] e^\alpha}
		{1+te^\alpha}(f_t)\cdot\alpha
\]
the limit of the last term at $x=0$ can be found using \eqref{limits} 
and leads to the constraint
\begin{equation}\label{linear}
	\iota(\xi_1)\left[hb + s H_V(f_{-1})\right] =
	\sum\nolimits_\beta \frac{\beta}{\beta(\xi_1)},
\end{equation}
summed over the roots $\beta$ of $\frz$. Its solutions are 
the critical points of the function of $\zeta\in\frt$
\[
\zeta \mapsto  \frac{1}{2}\left[hb(\zeta,\zeta) + s H_V(\zeta,\zeta)\right] -
		\textstyle\sum_\beta \log|\beta(\zeta)|.
\]  
This function is real-valued for $s=0$, blows up on the walls of each 
Weyl chamber of $\frz$ and is dominated by the quadratic term at large $\zeta$, 
so a minimum must exist inside the chamber. Further, the Hessian 
\[
hb + \sum\nolimits_\beta \frac{\beta^{\otimes 2}}{\beta(\zeta)^2}
\]
is positive-definite, so the minimum is non-degenerate and the $s$-perturbed
equation can also be solved for small $s$. Continuing the solution to 
higher order in $x$, we get an recursive family of equations in each 
degree $k>1$
\begin{equation}\label{recursion}
\iota(\xi_k)\left(hb + sH_V + \sum\nolimits_\beta \frac
	{\beta^{\otimes 2}}{\beta(\xi_1)^2}\right) = 
		(\text{expression in } \xi_j,\:\: j<k),
\end{equation}
which can be solved by reason of the same non-degeneracy.  
\end{proof}

\vskip.5cm
\small{
\noindent\textsc{C.\ Teleman:} DPMMS, CMS, Wilberforce Road, Cambridge 
CB2 1TP, UK. \texttt{teleman@dpmms.cam.ac.uk}\\
\textsc{C.T.\ Woodward:} Dept.~of Mathematics, Hill Center, Rutgers 
University, 110 Frelinghuysen Road, Piscataway NJ 08854, USA. \texttt{ctw@math.rutgers.edu}}
\end{document}